\documentclass[a4paper, 12pt]{article}

\usepackage[koi8-r,cp1251]{inputenc}
\usepackage{amsfonts,amssymb,amsmath, hyperref}
\usepackage[final]{epsfig}
\usepackage{graphicx}

\begin{document}

\begin{center}
\bf{{\Large Boundedness of  Hausdorff  operators on Hardy spaces over
 homogeneous   spaces of Lie groups }}
\end{center}

\begin{center}
A. R. Mirotin
\end{center}

\

\begin{center}
airotin@yandex.ru
\end{center}

\

Abstract. The aim of this note is to give the boundedness conditions for  Hausdorff operators on Hardy spaces $H^{1}$ with the norm defined via  $(1,q)$ atoms over homogeneous   spaces of Lie groups with doubling property  and to apply  results we obtain to
generalized Delsarte operators and to Hausdorff  operators over multidimensional spheres.

\

Key words. Hausdorff operator, Lie group, homogeneous   space, Hardy space, generalized  shift operator of Delsarte.

\

MSC classes: 43A85, 47G10, 22E30

\

\section{Introduction}

One-dimensional Hausdorff operators were introduced by Hardy \cite[Section 11.18]{H} as a transformations
of functions of a continuous variable analogous to the regular Hausdorff transformations for sequences and series. Although occasionally one-dimensional Hausdorff operators appeared before 2000   (see \cite{X1} and \cite{X2}),  the modern development of this theory begins with the work of Liflyand and  M\'{o}ricz \cite{LM} where   Hausdorff  operators on one-dimensional Hardy space were considered. The multidimensional case was studied in \cite{LL}. For more details of the  development of the theory of Hausdorff  operators up to 2014 see \cite{Ls}, and \cite{CFW}.

Hausdorff operators on the Hardy space $H^{1}$ over homogeneous spaces of locally compact groups were first introduced
by the author in \cite{JBSU} for the case of doubling measures, and in \cite{AnMath} for the case of locally  doubling measures. The case of locally compact groups was considered earlier in \cite{JMAA}.  The aim of this note is to improve and generalize results from \cite{JBSU} to the case of Hardy spaces $H^{1}(G/K)$ with the norm defined via  $(1,q)$ atoms when $G$ is a Lie group and to apply  results we obtain to
generalized Delsarte operators and to Hausdorff  operators over multidimensional spheres.

\section{The main result}

 Let $G$ be a locally compact metrizable group with left invariant  metric $\rho$ and left Haar  measure $\nu$. We assume that the following \textit{doubling condition} in a  sense of \cite{CW} holds:

there exists  a constant $C$ such that
$$
\nu(B(x,2r))\leq C \nu(B(x,r))
$$
for each $x\in G$ and $r > 0$;
 here $B(x,r)$ denotes the ball of radius $r$ around $x.$

 The \textit{doubling constant }  is the smallest constant
$C\geq  1$ for which the last inequality is  valid. We denote this constant by $C_\nu.$ Then  for each $x\in G, k\geq 1$ and $r > 0$
$$
\nu(B(x,kr))\leq C_\nu k^{d} \nu(B(x,r)),\eqno(D)
$$
where $d=\log_2C_\nu$ (see, e.g., \cite[p. 76]{HK}). The number $d$   takes
the role of a "dimension"  for a doubling metric measure space $G$.

 Homogeneous group in a sense of Folland  and Stein  \cite{FS}
(i.e., a connected simply connected Lie group  $G$ whose Lie algebra is equipped
with dilations) enjoys the doubling condition and $C_{\nu}=2^Q,$ where $Q$ stands for the homogeneous dimension of $G$ \cite[Lemma 3.2.12]{FR}.
A complete
Riemannian manifold, with Ricci curvature nonnegative outside a compact
subset of the manifold, satisfies the doubling condition, as demonstrated in \cite[Lemma 1.3]{LT}.
Compact Lie groups endowed with Riemann metric and Haar measure satisfy the doubling condition, too \cite[p. 588, Example (7)]{CW}.
For complete noncompact
manifolds with nonnegative Ricci curvature, the doubling property
for the volume measure follows from the volume comparison inequality of
Bishop and Gromov \cite[Theorem 10.6.6]{BBI}.

 We denote by $\mathrm{Aut}(G)$  the space  of all topological  automorphisms of $G$ endowed with its natural topology $\mathcal{T}_\beta$ \cite[Ch. X, \S 3, n  5]{Burb},  $\mathcal{L}(Y)$ denotes the space of linear bounded operators on a normed space $Y$.

 Let $K$ be a compact subgroup of  $G$ with normalized Haar measure $\beta$.  Consider the quotient space
$G/K$ of left cosets $\dot x:=xK=\pi_K(x)$ ($x\in G$) where $\pi_K:G\to G/K$ stands for a natural projection. We shall assume that the  measure $\nu$ is normalized in such a way that (generalized) Weil's formula
$$
\int_Gg(x)dx=\int_{G/K}\left(\int_Kg(xk)dk\right)d\lambda(\dot x) \eqno(1)
$$
holds for all $g\in L^1(G)$, where $\lambda$ denotes some left-$G$-invariant measure on $G/K$ (see
\cite[Chapter VII, §2, No. 5, Theorem 2 ]{Bourb} and especially  remark c) after this theorem or \cite[Proposition 10.4.12]{HN}).
Here $G$-left invariance of   $\lambda$ means that $\lambda(xE)=\lambda(E)$ for every Borel subset $E$ of $G/K$ and for every $x\in G$. This measure is unique up to constant multiplier.

   Henceforth we  write $dx$ instead of $d\nu(x)$ and $dk$ instead of $d\beta(k)$. We shall write also $d\dot x$ instead of $d\lambda(\dot x)$.

The function $g:G\to \mathbb{C}$ is called \textit{right-$K$-invariant}  if $g(xk)=g(x)$ for all $x\in G$, $k\in K$. For such a function we put $\dot g(\dot x):=g(x)$. This definition is correct and for $g\in L^1(G)$ formula (1) implies that
$$
\int_Gg(x)dx=\int_{G/K}\dot g(\dot x)d\dot x \eqno(2)
$$
 (recall that $\int_Kdk=1$).

 The map $g\mapsto\dot g$ is a bijection between the set of all right-$K$-invariant functions on $G$ (all right-$K$-invariant functions from $L^1(G)$) and the set of all functions on $G/K$ (respectively functions from $L^1(G/K,\lambda)$).

Let an automorphism $A\in \mathrm{Aut}(G)$ maps $K$ onto itself. Since
$$
A(\dot x):= A(xK)=\{A(x)A(k): k\in K\}=
A(x)K=\pi_K(A(x))
$$
 we get a homeomorphism   $\dot A:G/K\to G/K,$  $\dot A(\dot x):=\pi_K(A(x)).$ Then for every
 right-$K$-invariant function $g$ on $G$ we have  $\dot g(\dot A(\dot x))=g(A(x)).$

We put
 $$
 \mathrm{Aut}_K(G):=\{\dot A : A\in \mathrm{Aut}(G), A(K)=K\}.
 $$

A $\nu$-measurable function $a$ on $G$ is called an $(1,q)$-\textit{atom} ($q\in (1,\infty]$) if

(i) the support of $a$ is contained in a ball $B(x, r)$;

(ii) $\|a\|_\infty \le \frac{1}{\nu(B(x,r))}$ if $q=\infty$, and

$\|a\|_q \le \nu(B(x,r))^{\frac{1}{q}-1}$  if $q\in (1,\infty)$\footnote{As usual, $\|\cdot\|_q$ denotes the $L^q$ norm.};

(iii) $\int_G a(x)d\nu(x) = 0$.

In case $\nu(G)<\infty$ we shall assume $\nu(G)=1$; in this case  the constant function
having value $1$ is also considered to be an atom.

Hereafter   by atom we mean an $(1,q)$-atom on $G$.

\textbf{Definition 1.} \cite{JBSU}, \cite{AnMath}. We define the \textit{Hardy space} $H^{1}(G/K)=H^{1,q}(G/K)$\footnote{It is known that $H^{1,q}(G/K)$ does not depend on $q\in (1,\infty]$ \cite[Theorem A, p. 592]{CW}. We write $H^{1,q}(G/K)$ instead of $H^{1}(G/K)$ in order to stress the fact that we use the  norm $\|\cdot\|_{H^{1,q}(G/K)}$ described below.} as a space of such functions $f$  on  $G/K$ that  admit an atomic decomposition  of the form
$$
f=\sum_{j=1}^\infty \lambda_j\dot a_j
$$
where $a_j$ are  right-$K$-invariant  $(1,q)$-atoms on $G$ and $\sum_{j=1}^\infty |\lambda_j|<\infty$.
In this case,
$$
\|f\|_{H^{1,q}(G/K)}:=\inf\sum_{j=1}^\infty |\lambda_j|,
$$
and infimum is taken over all  decompositions above of $f$.

In other words, $f=\dot g$ where
$g=\sum_{j=1}^\infty \lambda_ja_j$,
 $a_j$ are  right-$K$-invariant  $(1,q)$-atoms on $G$, and $\sum_{j=1}^\infty |\lambda_j|<\infty$.
Moreover,
  $\|f\|_{H^{1,q}(G/K)}=\|g\|_{H^{1,q}(G)}$.

\textbf{Remark 1}. Real Hardy spaces over compact connected (not necessary quasi-metric) Abelian groups  were  defined in \cite{Indag}.

\textbf{Proposition 1}. \cite{AnMath}. \textit{Let   $G\ne K$. Then the space $H^{1,q}(G/K)$ is nontrivial and  Banach.}

\textbf{Definition 2.} \cite{JBSU}. Let  $(\Omega,\mu)$ be  a measure space,  $(\dot A(u))_{u\in \Omega}\subset \mathrm{Aut}_K(G)$  a family of homeomorphisms of $G/K$,  and
$\Phi$ a measurable function on $(\Omega,\mu)$.  For   a  Borel measurable function $f$ on $G/K$ we define a  \textit{Hausdorff operator on} $G/K$ as follows
$$
(\mathcal{H}_{\Phi,\dot A} f)(\dot x):=\int_{\Omega} \Phi(u)f(\dot A(u)(\dot x))d\mu(u).
$$

 For the proof of our main result the next two lemmas  are crucial.

\textbf{Lemma 1.}  \cite{JMAA}. \textit{ Let $(\Omega,\mu)$ be  $\sigma$-compact quasi-metric space
with   positive Radon measure $\mu,$  $(X,m)$ be a measure space and $\mathcal{F}(X)$
be some Banach  space of  $m$-measurable functions on $X$. Assume  that the convergence
of a sequence strongly in  $\mathcal{F}(X)$ yields the convergence of some subsequence
to the same function for $m$-almost all $x\in X$. Let $F(u,x)$ be a function such that
$F(u,\cdot)\in \mathcal{F}(X)$ for $\mu$-almost all $u\in \Omega$ and the
map $u\mapsto F(u,\cdot):\Omega\to \mathcal{F}(X)$ is Bochner integrable with respect
to $\mu$. Then for $m$-almost all $x\in X$}
$$
\left((B)\int_\Omega F(u,\cdot)d\mu(u)\right)(x)=\int_\Omega F(u,x)d\mu(u).
$$

\textbf{Lemma 2.} \textit{Let $G$ be a (finite dimensional real or complex) connected
 Lie group with left invariant Riemann metric $\rho$. Then every automorphism  $\varphi\in  \mathrm{Aut}(G)$  is Lipschitz with Lipschitz constant $\|(d\varphi)_e\|$.}

Proof.  Let $T_a(G)$ denotes the tangent space for $G$ at the point $a\in G$. Let $L_a:x\mapsto ax$ be the left translation in $G$. Then the tangent map $l_a:=(dL_a)_e:T_e(G)\to T_a(G)$ is a bijection. We fix  the Euclidean norm $\|\cdot\|$  in $T_e(G)$ and  introduce the norm in $T_a(G)$ by the rule  $\|X_a\|:= \|X_e\|$ if $X_a=l_a(X_e)$, $X_e\in T_e(G)$, $a\in G$.

 As is well known, for every $p, q\in G$
$$
\rho(p, q)=\inf\limits_\alpha\int_0^1\left\|\alpha'(t)\right\|dt
$$
where infimum is taken over all  piecewise smooth curves  $\alpha$ from $[0,1]$ to $G$ with  $\alpha(0)=p,$  $\alpha(1)=q$ ($\alpha'(t)$ stands, as usual, for the tangent vector to $\alpha$ at the point $\alpha(t)$).
Since  $\varphi\in  \mathrm{Aut}(G)$, the formula $\beta=\varphi\circ\alpha$ gives the general form of all piecewise smooth  curves in $G$ with $\beta(0)=\varphi(p)$ and $\beta(1)=\varphi(q)$. Thus, by the chain rule
$$
\rho(\varphi(p), \varphi(q))=\inf\limits_\alpha\int_0^1\left\|(\varphi\circ\alpha)'(t)\right\|dt
$$
$$=
\inf\limits_\alpha\int_0^1\left\|(d\varphi)_{\alpha(t)}\alpha'(t)\right\|dt
\le \inf\limits_\alpha\int_0^1\left\|(d\varphi)_{\alpha(t)}\|\|\alpha'(t)\right\|dt.
$$
It is known (see, e.g., \cite{BourbLie}) that for every left invariant vector field $X$ on $G$
 (this means that $X_a=l_a(X_e)$ for all $a\in G$) the vector field $(d\varphi)(X)$\footnote{In \cite{BourbLie} the map $d\varphi$ is denoted by $L(\varphi)$} is  left invariant, too. In other wards,  $(d\varphi)_a(X_a)=l_a(d\varphi)_e(X_e)$, i.e., $(d\varphi)_a=l_a((d\varphi)_e)l_a^{-1}$ and therefore $\|(d\varphi)_a\|=\|(d\varphi)_e\|$ for all $a\in G$. The result follows.

Now we are in a position to prove the next

\textbf{Theorem 1.} \textit{Let $G$ be a (finite dimensional real or complex) connected
 Lie group with left invariant Riemann metric $\rho$ and left Haar measure $\nu$  such that    the space $(G,\rho,\nu)$ is doubling. Let $(\Omega,\mu)$ be  $\sigma$-compact quasi-metric space
with   positive Radon measure $\mu$, and let $q\in (1,\infty]$. If
$$
\|\Phi\|_{A,q}:=\int_\Omega |\Phi(u)|(\mathrm{mod}A(u))^{-\frac{1}{q}}k(u)^{(1-\frac{1}{q})d}d\mu(u)<\infty
$$
where   $k(u):=\|(d(A(u)^{-1}))_e\|$, then the operator $\mathcal{H}_{\Phi,\dot A}$ is bounded on the  space $H^{1,q}(G/K)$ and }
$$
\|\mathcal{H}_{\Phi,\dot A}\|_{\mathcal{L}(H^{1,q}(G/K))}\leq C_\nu^{1-\frac{1}{q}}\|\Phi\|_{A,q}.
$$

Proof.  If we set $X=G/K$  and $m=\lambda$
 the pair $(X,m)$ satisfies  the conditions of Lemma 1 with $H^{1,q}(G/K)$ in place of $\mathcal{F}(X)$. Indeed, let $f_n=\dot g_n\in H^{1,q}(G/K)$, $f=\dot g\in H^{1,q}(G/K)$, and $\|f_n-f\|_{H^{1,q}(G/K)}\to 0$ ($n\to\infty$). Since
$$
\|f_n-f\|_{L^1(G/K)}=\int_{G/K}|\pi_K(g_n- g)|d\lambda
$$
$$
=\int_{G}|g_n(x)-g(x)|dx\leq \|g_n-g\|_{H^{1,q}(G)}=\|f_n-f\|_{H^{1,q}(G/K)}\to 0
$$
(by H\"{o}lder inequality $\|a\|_1\le 1$ for each atom $a$), there is a subsequence of $f_{n}$ that converges to $f$ $\lambda$-a.e.

Then Definition 2 and Lemma 1 imply for $f\in H^{1,q}(G/K)$ that
$$
\mathcal{H}_{\Phi, \dot A}f= \int_\Omega \Phi(u) f\circ \dot A(u)d\mu(u),
$$
the Bochner integral (recall that $ H^{1,q}(G/K)$ is a subspace of  $L^{1}(G/K,\lambda)$ \cite[p. 592]{CW}, and thus we identify functions that equal $\lambda$-a.e.).

Therefore (below $f=\dot g$)
$$
\|\mathcal{H}_{\Phi,\dot A}f\|_{H^{1,q}(G/K)}\leq \int_\Omega |\Phi(u)|\|f\circ \dot A(u)\|_{H^{1,q}(G/K)}d\mu(u)
$$
$$
=\int_\Omega |\Phi(u)|\|g\circ A(u)\|_{H^{1,q}(G)}d\mu(u).
$$
If $g=\sum_{j=1}^\infty \lambda_ja_j$  then
$$
g\circ A(u)=\sum_{j=1}^\infty\lambda_ja_j\circ A(u).\eqno(3)
$$
We claim that
$$
b_{j,u}:=C_\nu^{\frac{1}{q}-1}(\mathrm{mod}(A(u))^{\frac{1}{q}}k(u)^{(\frac{1}{q}-1)s}a_j\circ A(u)
$$
 is an atom, too. Indeed, Lemma 2 implies that
 $$
A(u)^{-1}(B(x,r))\subseteq B(x',k(u)r),
 $$
 where $x'=A(u)^{-1}(x)$. If $a_j$ is supported in $B(x_j,r_j)$ then $b_{j,u}$ is supported in $B(x'_j,k(u)r_j)$. So the condition (i) holds for  $b_{j,u}$.

 Next, by  the  property (D) we have
 $$
\nu(B(x_j,k(u)r_j))\le C_\nu k(u)^d\nu(B(x_j,r_j)).
 $$
This estimate yields in view of (ii) and the left invariance of $\rho$ and $\nu$ that
 $$
 \|a_j\circ A(u)\|_q=\left(\int_G|a_j\circ A(u)|d\nu\right)^\frac{1}{q}=(\mathrm{mod}(A(u))^{-\frac{1}{q}}\|a_j\|_q
 $$
 $$
 \le (\mathrm{mod}(A(u))^{-\frac{1}{q}}(\nu(B(x_j,r_j)))^{\frac{1}{q}-1}
  \le (\mathrm{mod}(A(u))^{-\frac{1}{q}}\left(\frac{\nu(B(x_j,k(u)r_j))}{C_\nu k(u)^d}\right)^{\frac{1}{q}-1}
  $$
  $$
  =\left(C_\nu^{\frac{1}{q}-1}(\mathrm{mod}(A(u))^{\frac{1}{q}}k(u)^{(\frac{1}{q}-1)d}\right)^{-1} (\nu(B(x_j',k(u)r_j)))^{\frac{1}{q}-1}.
 $$
 Thus, condition (ii) holds for $b_{j,u}$, too. Finally, the validity of (iii) follows from
 \cite[VII.1.4, formula (31)]{Bourb}.

 Since formula (3) can be rewritten in the form
 $$
g\circ A(u)=\sum_{j=1}^\infty\left(\lambda_jC_\nu^{1-\frac{1}{q}}(\mathrm{mod}(A(u))^{-\frac{1}{q}}k(u)^{(1-\frac{1}{q})d}\right)b_{j,u},
$$
 we have
$$
\|g\circ A(u)\|_{H^{1,q}(G)}\leq C_\nu^{1-\frac{1}{q}}(\mathrm{mod}(A(u))^{-\frac{1}{q}}k(u)^{(1-\frac{1}{q})d}\sum_{j=1}^\infty|\lambda_j|.
$$
It follows that (recall that $f=\dot g$)
$$
\|g\circ A(u)\|_{H^{1,q}(G)}\leq C_\nu^{1-\frac{1}{q}}(\mathrm{mod}(A(u))^{-\frac{1}{q}}k(u)^{(1-\frac{1}{q})d}\|g\|_{H^{1,q}(G)}
$$
$$
=
C_\nu^{1-\frac{1}{q}}(\mathrm{mod}(A(u))^{-\frac{1}{q}}k(u)^{(1-\frac{1}{q})d}\|f\|_{H^{1,q}(G/K)}.
$$

Therefore
$$
\|\mathcal{H}_{\Phi,\dot A}\|_{\mathcal{L}(H^{1,q}(G/K))}\leq C_\nu^{1-\frac{1}{q}}\int_\Omega |\Phi(u)|(\mathrm{mod}A(u))^{-\frac{1}{q}}k(u)^{d(1-\frac{1}{q})}d\mu(u)
$$
and the proof is complete.

Setting in Theorem 1 $\Omega=\mathbb{Z}_+$ with counting measure $\mu$ we have the next result
for discrete Hausdorff operators.

\textbf{Corollary 1.} \textit{Let $(G,\rho,\nu)$ and $K$ be as in  the  Theorem 1, $(\dot A(n))_{n\in \mathbb{Z}_+}\subset \mathrm{Aut}_K(G)$, and  $q\in (1,\infty]$. If $\Phi:\mathbb{Z}_+\to \mathbb{C}$ be such that
 $$
\|\Phi\|_{A,q}:=\sum\limits_{n=0}^\infty |\Phi(n)|(\mathrm{mod}A(n))^{-\frac{1}{q}}k(n)^{(1-\frac{1}{q})d}<\infty,
$$
 then the $\mathrm{discrete\ Hausdorff\ operator}$
 $$
 \mathcal{H}_{\Phi,\dot A}f(\dot x):=\sum\limits_{n=0}^\infty \Phi(n)f(\dot A(n)(\dot x))
 $$ is bounded on the  space $H^{1,q}(G/K)$ and }
$$
\|\mathcal{H}_{\Phi,\dot A}\|_{\mathcal{L}(H^{1,q}(G/K))}\leq C_\nu^{1-\frac{1}{q}}\|\Phi\|_{A,q}.
$$

As a special case of Theorem 1 for $K=\{e\}$ ($e$ denotes the unit of $G$) one has the

\textbf{Corollary 2.} \textit{ Let Let $(G,\rho,\nu)$ and $(\Omega,\mu)$   be as in  the  Theorem 1, $(A(u))_{u\in \Omega}\subset \mathrm{Aut}(G)$, and $q\in (1,\infty]$. If $\|\Phi\|_{A,q}<\infty$
then the operator $\mathcal{H}_{\Phi, A}$ is bounded on  $H^{1,q}(G)$ and }
$$
\|\mathcal{H}_{\Phi, A}\|_{\mathcal{L}(H^{1,q}(G))}\leq C_\nu^{1-\frac{1}{q}}\|\Phi\|_{A,q}.
$$

\textbf{Remark 2}. The condition $\|\Phi\|_{A,q}<\infty$ is not necessary for boundedness of
$\mathcal{H}_{\Phi, A}$ in Hardy space as the following simple example shows\footnote{The sufficient boundedness conditions from \cite {LL} and \cite{CZ} are also not met in this example.}.  Consider the Hausdorff operator
$$
(\mathcal{H}_1f)(x):=\int_{\Omega}f(u_1x_1,\dots,u_nx_n)du
$$
in $ H^1(\mathbb{R}^n)$.  Here $G=\mathbb{R}^n$, $\Omega=\{u\in\mathbb{R}^n: u_j\ne 0 \mbox{  for } j=1,\dots,n \}$, $\mu$ and $\nu$ are  Lebesgue measures on $\Omega$   and $\mathbb{R}^n$ respectively, $K=\{0\}$, $A(u)(x)=A_ux$, where $A_u=\mathrm{diag}\{u_1,\dots,u_n\}$ ($x\in \mathbb{R}^n$ a column vector, $u\in \Omega$), $\Phi=1$, $d=n$. The necessary moment condition  $\int_{\mathbb{R}^n}f(u)du=0$ for functions from $H^1(\mathbb{R}^n)$  yields that $\mathcal{H}_1f=0$ for all $f\in H^1(\mathbb{R}^n)$. On the other hand, here $\mathrm{mod}A(u)=|\det A_u|=|u_1\dots u_n|$ \cite[Subsection VII.1.10, Corollary 1]{Burb1}, $(dA(u)^{-1})_0X=A_u^{-1}X$ ($X\in \mathbb{R}^n$), $k(u)=\|A_u^{-1}\|=(\sum_{j=1}^nu_j^{-2})^{1/2}\ge n^{1/2}|u_1\dots u_n|^{-1/n}$. Then
$$
\|\Phi\|_{A,q}=\int_{\Omega}(\mathrm{mod}A(u))^{-\frac{1}{q}}k(u)^{(1-\frac{1}{q})d}du\ge n^{\frac{1}{2}(1-\frac{1}{q})n}\int_{\Omega}\frac{du}{|u_1\dots u_n|}
=\infty.
 $$

\section{Examples}

\subsection{Generalized   shift operator of Delsarte}

Let $G$ be as above and $\mathfrak{A}$  a compact subgroup of $\mathrm{Aut}(G)$ with normalized Haar measure $m$.
  Recall that the generalized   shift operator of Delsarte  \cite{Delsarte}, \cite[Ch. I, \S 2]{Lev} (also  the  terms “generalized  translation operator of Delsarte”, or “generalized displacement operator of Delsarte” are used) is defined to be
$$
T^xf(h)=\int_{\mathfrak{A}} f(ha(x))dm(a)\quad (x,h\in G).
$$
Since the group $G$ acts on $G/K$, one can define a generalization of this   operator   to $G/K$ as follows.
   Let $\Omega:=\{u\in \mathfrak{A}: u(K)=K\}$. Then $\Omega$ is a compact subgroup of $\mathfrak{A}$. We denote by $\mu$ the normalized Haar measure of $\Omega$ and put for a Borel measurable function $f$ on $G/K$
  $$
  T^{\dot x}f(h):=\int_\Omega f(h\dot u(\dot x))d\mu(u)\quad (\dot x\in G/K,h\in G).
  $$
  Let $h$ be fixed and $L^hf(\dot x):=T^{\dot x}f(h)$.
Then $L^h=\mathcal{H}_1\tau_h$, where
$$
\mathcal{H}_1f(\dot x):=\int_\Omega f(\dot u(\dot x))d\mu(u)
$$
is a Hausdorff operator on $G/K$ with $\Phi(u)=1$ and $A(u)=u$, and $\tau_h f(\dot x):=f(h\dot x)$.  Note that  $\mathrm{mod}$ is a continuous homomorphism from $\mathrm{Aut}(G)$ to the multiplicative group $(0,\infty)$. Since $\Omega$ is a compact group, it follows that  $\mathrm{mod}(\Omega)=\{1\}$.  Assume that
the doubling conditions for the Lie group $G$ holds and   $\Omega$ is  quasi-metric. Then the operator $\mathcal{H}_1$ is bounded on $H^{1,q}(G/K)$  by Theorem 1 and
$$
\|\mathcal{H}_{1}\|\le C_\nu^{1-\frac{1}{q}}\int_\Omega k(u)^{(1-\frac{1}{q})s}d\mu(u)
$$
where $k(u)=\|(d(u^{-1}))_e\|$. Since $\tau_h$ is an isometry of $H^{1,q}(G/K)$, we conclude that  the operator $L^h$ is bounded on $H^{1,q}(G/K)$  and
$$
\|L^h\|\le C_\nu^{1-\frac{1}{q}}\int_\Omega k(u)^{(1-\frac{1}{q})d}d\mu(u).
$$

\subsection{Hausdorff operators on the unit sphere in $\mathbb{R}^n$}

Consider the unit sphere $\mathbb{S}^{n-1}\subset \mathbb{R}^n$ (the case $n=3$ was considered in \cite{AnMath}).

The compact group $G=SO(n)$ acts on  $\mathbb{S}^{n-1}$ transitively by restriction of the natural action of $GL(n,\mathbb{R})$ on $\mathbb{R}^n$. It is known that the
isotropy subgroup $K$ of the point $e_n:=(0,\dots,1)\in \mathbb{S}^{n-1}$ consists of all elements in $SO(n)$ of the form
$$
\widetilde{a}:=\left(\begin{array}{cc}
 a&\mathbf{0}^\top\\
\mathbf{0}&1
\end{array}\right),
$$
where $\mathbf{0}=(0,\dots,0)\in \mathbb{R}^{n-1}$, $a\in  SO(n-1)$. Hence we identify
$\mathbb{S}^{n-1}$ with the homogeneous space $SO(n)/K$. Let $s\in \mathbb{S}^{n-1}$. If a matrix $x(s)\in SO({n-1})$ satisfies $s=x(s)e_n^\top$ we can identify the point $s$ with the coset $\dot x(s):=x(s)K$.

 Consider the set of automorphisms of  $G=SO(n)$ of the form
 $$
 A(u)(x)=\widetilde{u}^{-1}x\widetilde{u},\quad u\in O(n-1).
 $$
 Since every mapping $x\mapsto u^{-1}xu$ with $u\in O(n-1)$  maps $SO(n-1)$ onto itself (being a connected component of unit in $O(n-1)$ the group $SO(n-1)$ is a normal subgroup of $O(n-1)$), we have in our case that all automorphisms  $A(u)$ where $u\in O(n-1)$
map $K$ onto $K$. Then by definition the coset
$$
\dot A(u)(\dot x(s))= \pi_K(\widetilde{u}^{-1}x(s)\widetilde{u})
$$
can be identified with the  point $$
\widetilde{u}^{-1}x(s)\widetilde{u}e_n^\top=\widetilde{u}^{-1}x(s)e_n^\top=\widetilde{u}^{-1}s=(u^{-1}s',s_n)
$$
($s':=(s_1,\dots,s_{n-1})$) of $\mathbb{S}^{n-1}$.

Thus, Definition 2 takes the form (we put $x=x(s)$ in this definition and  identify the coset $\dot x(s)$ with a column vector $s\in \mathbb{S}^{n-1}$)
$$
(\mathcal{H}_{\Phi, \mu}f)(s) =\int_{O(n-1)} \Phi(u)f(u^{-1}s',s_n)d\mu(u)\eqno(4)
$$
where  $\mu$ stands for a (regular Borel) measure on $O(n-1)$ and $f$ is a  Borel measurable   function on $\mathbb{S}^{n-1}$.

Note that the point $(u^{-1}s',s_n)$ runs over the cross-section of  $\mathbb{S}^{n-1}$ by the hyperplane $\{x=s_n\}\subset \mathbb{R}^n$ (which contains $s$) orthogonal to the last coordinate axis  when $u$ runs over $O(n-1)$. So (4) looks as  a "horizontal slice transform" on $\mathbb{S}^{n-1}$ and the function $\mathcal{H}_{\Phi, \mu}f$ depends on $s_n\in[-1,1]$ only.

To apply Theorem 1 first we shaw that $k(u)=1$  for $u\in O(n-1)$. Indeed, $k(u)=\|d(A(u^{-1}))_{1_n}\|$ (here $1_n$ stands for the unit $n\times n$ matrix). It is easy to verify that for every $X\in \mathfrak{so}(n)$, the Lie algebra of $SO(n)$
$$
d(A(u^{-1}))_{1_n}X=\widetilde{u}X\widetilde{u}^{-1}.
$$
On the other hand, $\widetilde{u}\in O(n)$ for $u\in O(n-1)$. Thus,
$$
\|d(A(u^{-1}))_{1_n}\|=
\max\limits_{\|X\|=1,\|Y\|=1}|\langle d(A(u^{-1}))_{1_n}X,d(A(u^{-1}))_{1_n}Y\rangle|=
$$
$$
\max\limits_{\|X\|=1,\|Y\|=1}|\langle\widetilde{u}X\widetilde{u}^{-1},\widetilde{u}Y\widetilde{u}^{-1}\rangle|
=\max\limits_{\|X\|=1,\|Y\|=1}|\langle X,Y\rangle|=1
$$
(here $\langle\cdot,\cdot\rangle$ stands for the Euclidean inner product).
 Since $SO(n)$ is compact, it is doubling \cite{CW}.  Next, since $SO(n)$ is unimodular, we  get that $\mathrm{mod} A=1$ for all $A\in \mathrm{Aut}(SO(n))$. So if $\Phi\in L^1(O(n-1),\mu)$ \cite[Theorem 1]{JBSU} yields that the operator (4) is bounded on  $L^p(\mathbb{S}^{n-1})$   and $\|\mathcal{H}_{\Phi,\mu}\|_{\mathcal{L}(L^p(\mathbb{S}^{n-1}})\le \|\Phi\|_{L^1(\mu)}$. Moreover,   Theorem 1 yields that
$$
\|\mathcal{H}_{\Phi,\mu}\|_{\mathcal{L}(H^{1,q}(\mathbb{S}^{n-1}))}\le C_\nu^{1-\frac{1}{q}}\|\Phi\|_{L^1(\mu)}
$$
where  $C_\nu$ is the doubling constant for  $SO(n)$.

In closing let us consider the following special case. Let $\Phi=1$ and $m$ be a Haar measure of the (compact) group $O(n-1)$. Then for every $f\in H^1(\mathbb{S}^{n-1})$ the function
$$
(\mathcal{H}_{1,m}f)(s)=\int_{O(n-1)}f(u^{-1}s',s_n)dm(u)
$$
belongs to $H^1(\mathbb{S}^{n-1})$. On the other hand, this function  depends on $s_n$ only.
Indeed,  $(s',s_n)\in \mathbb{S}^{n-1}$ if and only if  $s'$ belongs to the sphere  $\mathbb{S}^{n-2}_r$
centered at $0\in\mathbb{ R}^{n-1}$  of radius $r=\sqrt{1-s_n^2}$. Fix  $s'_0\in\mathbb{S}^{n-2}_r$. Since $SO(n-1)$ acts transitively on $\mathbb{S}^{n-2}_r$, for every
$s'\in\mathbb{S}^{n-2}_r$ there is such $v\in SO(n-1)$ that $vs'_0=s'$. Taking into account that $O(n-1)$ is unimodular, we get
$$
(\mathcal{H}_{1,m}f)(s)=\int_{O(n-1)}f(uvs'_0,s_n)dm(u)=\int_{O(n-1)}f(us'_0,s_n)dm(u)
$$
which completes the proof.

\textbf{Acknowledgments}.   This work was supported by the State Program of Scientific Research of the Republic of Belarus.


\begin{thebibliography}{99}

\bibitem{H}
 G. H. Hardy, \textit{Divergent Series}. Clarendon Press (Oxford, 1949).


\bibitem{X1}
Richard R. Goldberg,  Certain operators and Fourier transforms on $L^2$. \textit{Proc. Amer. Math. Soc.} \textbf{10} (1959), 385--390.

 \bibitem{X2}
 C. Georgakis, The Hausdorff mean of a Fourier-Stieltjes transform. \textit{Proc. Amer. Math. Soc.} \textbf{116} (1992), 465--471.



\bibitem{LM}
 E. Liflyand,  F. M\'{o}ricz, The Hausdorff operator is bounded on the
real Hardy space $H^1(\mathbb{R})$.\textit{ Proc. Amer. Math. Soc.}  \textbf{128} (2000) 1391--1396. DOI: https://doi.org/10.1090/S0002-9939-99-05159-X.



\bibitem{LL}
 A. Lerner and  E. Liflyand, Multidimensional Hausdorff operators on the
real Hardy space, \textit{J. Austr. Math. Soc.}, \textbf{83}  (2007), 79--86.





\bibitem{Ls}
 E. Liflyand, Hausdorff operators on Hardy spaces. \textit{Eurasian Math. J.}, \textbf{4}(4)
 (2013), 101--141.

\bibitem{CFW}
J. Chen, D. Fan, and S. Wang, Hausdorff operators on Euclidean space (a survey article). \textit{Appl. Math. J. Chinese Univ. Ser. B.}, \textbf{28}(4)  (2013), 548-564. DOI: https://doi.org/10.1007/s11766-013-3228-1.





\bibitem{JBSU}
A. R. Mirotin, Hausdorff operators on homogeneous spaces of locally
compact groups, Journal of the Belarusian State University. Math. Mech., No 2 (2020), 28--35.
 Letter to the editors  \textit{Ibid.}, No 3, p. 92 (2020)


\bibitem{AnMath}
A. R. Mirotin, Hausdorff  operators
on real Hardy spaces $H^1$ over homogeneous  spaces with local doubling property, Anal. Math. (to appear).


\bibitem{JMAA}
 A. R. Mirotin, Boundedness of Hausdorff operators on Hardy spaces $H^1$
over locally compact groups, \textit{J. Math. Anal. Appl.},  \textbf{473} (2019), 519 -- 533.
 DOI 10.1016/j.jmaa.2018.12.065. Preprint arXiv:1808.08257v2
[math.FA] 1 Sep 2018.




\bibitem{CW}
 R. R. Coifman, G. Weiss, Extensions of Hardy spaces and their use
in analysis, \textit{Bull. Amer. Math. Soc}. \textbf{83}   (1977), 569 -- 645.

\bibitem{HK}
Heinonen J.,  Koskela P., Shanmugalingam N.,  Tyson JT. \textit{Sobolev Spaces on Metric Measure Spaces. An Approach Based on Upper Gradients}. Cambridge: Cambridge University Press; 2015. 434 p. DOI:https://doi.org/10.1017/CBO9781316135914.



\bibitem{FS}
G.B. Folland,  E.M. Stein, Hardy Spaces on Homogeneous Groups. Mathematical
Notes, vol. 28. Princeton University Press, Princeton (1982)


\bibitem{FR}
V. Fischer, M. Ruzhansky,  Quantization on nilpotent Lie groups. Progress in
Mathematics, vol. 314. Birkhauser, Basel (2016).



\bibitem{LT}
Li, P. and Tam, L.-F.  Green’s functions, harmonic functions, and volume
comparison. \textit{J. Differential Geom.}, \textbf{41}(2) (1995), 277--318.

\bibitem{BBI}
Burago, D., Burago, Yu., and Ivanov, S. A Course in Metric Geometry.
Graduate Studies in Mathematics, vol. 33. Providence, RI: American Mathematical
Society  (2001).

\bibitem{Burb}
N. Bourbaki, General topology. Chapters 5--10,  Berlin Heidelberg, Springer-Verlag, 1989.


\bibitem{Bourb}
N. Bourbaki, Elements de mathematique. 2nd edn. Integration, Chaps.
I-IV. Hermann,  Paris, 1965.

\bibitem{HN}
J. Hilgert, K.-H. Neeb, Structure and Geometry of Lie Groups, Springer, New York - Dordrecht - Heidelberg - London (2012)






\bibitem{Indag}
 A. R. Mirotin,  On the general form of linear functionals on
the Hardy spaces $H^1$ over compact Abelian groups and
some of its applications. \textit{Indag. Math}.  \textbf{28} (2017), 451 -- 462.
DOI:https://doi.org/10.1016/j.indag.2016.11.023.


















\bibitem{BourbLie}
N. Bourbaki, Lie Groups and Lie Algebras. Chapters 1-3,  Berlin Heidelberg, Springer-Verlag (1989).

\bibitem{CZ}
J. Chen, X. Zhu, Boundedness of multidimensional Hausdorff operators
on $H^1(\mathbb{R}^n)$, J. Math. Anal. Appl. 409 (2014), 428--434.


\bibitem{Burb1}
N. Bourbaki,  Elements de mathematique. Livre VI. Integration. 2nd ed., Ch. 1 -- 9.   Hermann, Paris (1965 -- 1969).


\bibitem{Delsarte}
J. Delsarte, Hypergroupes et operateurs de permutation et de transmutation, \textit{Colloques Internat. Centre Nat. Rech. Sci.}, \textbf{71}   (1956), 29--45.

\bibitem{Lev}
B.M. Levitan, The theory of generalized displacement operators, Moscow, Nauka (1973) (In Russian). English translation: B.~M.~ Levitan,
Generalized translation operators and some of their applications,
Jerusalem, Israel Program for Scientific Translations (1964).








\end{thebibliography}
\end{document}